\documentclass{article}
\usepackage{amsfonts}
\usepackage{amsmath}

\newtheorem{theorem}{Theorem}

\newtheorem{proposition}[theorem]{Proposition}

\input{tcilatex}
\begin{document}

\title{The structure constants of geometric structures on a local Lie group}
\author{Erc\"{u}ment H. Orta\c{c}gil}
\maketitle

\begin{abstract}
We define the structure constants of the canonical almost complex, almost
symplectic and Riemannian structures on a local Lie group.
\end{abstract}

\section{A short review of absolute parallelizm}

We will recall some formulas from Part I of [O1] that will be needed below
and refer to [O1], [O2] for further details.

Let $M$ be a smooth manifold with the principal frame bundle $%
F(M)\rightarrow M.$ Elements of $F(M)$ are $1$-jets of local diffeomorphisms
(called $1$-arrows in [O1]) with source at the origin of $o$ of $\mathbb{R}%
^{n},$ $n=\dim M,$ and target at some $x\in M.$ $M$ is parallelizable if and
only if the structure group $GL(n,\mathbb{R)}$ of $F(M)\rightarrow M$ can be
reduced to the identity. Such a reduction $w,$ which we fix once and for
all, is called the structure object in [O1] and called a framing by
topologists. For each $x\in M,$ $w$ assigns a unique $1$-arrow from $o$ to $%
x.$ In terms of coordinates around $x,$ $w$ is of the form $w_{(j)}^{i}(x).$
The paranthesis around the index $j$ indicates that $w_{(j)}^{i}(x)$ is not
subject to any transformation upon a coordinate change $(x)\rightarrow (y)$
on $M,$ since this index refers to the standard coordinates in $\mathbb{R}%
^{n}.$ We define $z(x)=w(x)^{-1},$ which is a $1$-arrow from $x\in M$ to $%
o\in \mathbb{R}^{n},$ so that

\begin{equation}
w_{(a)}^{i}(x)z_{j}^{(a)}(x)=\delta _{j}^{i}\text{ \ \ \ \ }%
z_{a}^{(i)}(x)w_{(j)}^{a}(x)=\delta _{(j)}^{(i)}
\end{equation}%
with the standard summation convention. Note that the linear map $\delta
_{j}^{i}$ is the identity map at the tangent space $T_{x}(M)$ whereas $%
\delta _{(j)}^{(i)}$ is the identity map at $T_{o}(\mathbb{R}^{n})=\mathbb{R}%
^{n}.$ The structure object $w$ induces a transitive groupoid $\varepsilon $
on $M$ whose $1$-arrows on $M$ are given by

\begin{equation}
\varepsilon _{j}^{i}(x,y)\overset{def}{=}w_{(a)}^{i}(y)z_{j}^{(a)}(x)
\end{equation}

We define

\begin{eqnarray}
\Gamma _{jk}^{i}(x)\overset{def}{=}\left[ \frac{\partial \varepsilon
_{j}^{i}(x,y)}{\partial y^{k}}\right] _{y=x} &=&z_{j}^{(a)}(x)\frac{\partial
w_{(a)}^{i}(x)}{\partial x^{k}}= \\
-\left[ \frac{\partial \varepsilon _{j}^{i}(x,y)}{\partial x^{k}}\right]
_{y=x} &=&-\frac{\partial z_{j}^{(a)}(x)}{\partial x^{k}}w_{(a)}^{i}(x) 
\notag
\end{eqnarray}

\begin{equation}
T_{jk}^{i}(x)\overset{def}{=}\Gamma _{jk}^{i}(x)-\Gamma _{kj}^{i}(x)
\end{equation}

Now $T=(T_{jk}^{i})$ is called the integrability object of $T$ in [O2]. $%
\Gamma =(\Gamma _{jk}^{i})$ defines a first order linear operator $%
\widetilde{\nabla }$ which may be interpreted as linear connection on the
tangent bundle of $M.$ With this interpretation, $T$ turns out to be the
torsion of this connection and is called the torsion of $(M,w)$ in [O1], a
misleading name modified in [O2]. The linear curvature $\mathfrak{R}$ is
defined by

\begin{equation}
\mathfrak{R}_{jk,r}^{i}\overset{def}{=}\widetilde{\nabla }_{r}T_{jk}^{i}
\end{equation}%
so that $\mathfrak{R}=0$ if and only if $T$ is $\widetilde{\nabla }$%
-parallel or equivalently $\varepsilon $-invariant by Proposition 5.5, [O1].
In this case, the $1$-arrows of the groupoid $\Upsilon $ defined by (2)
integrate to a transitive pseudogroup $\mathcal{G}$ on $M$ whose local
diffeomorphisms $y=f(x)$ are solutions of the nonlinear PDE

\begin{equation}
\frac{\partial y^{i}}{\partial x^{j}}=\varepsilon
_{j}^{i}(x,y)=w_{(a)}^{i}(y)z_{j}^{(a)}(x)
\end{equation}

Now we assume $\mathfrak{R}=0$ so that $T$ is $\varepsilon $-invariant. This
means

\begin{equation}
\varepsilon _{a}^{i}(x,y)T_{bc}^{a}(x)\varepsilon _{j}^{b}(y,x)\varepsilon
_{k}^{c}(y,x)=T_{jk}^{i}(y)
\end{equation}%
Substituting (2) into (7) gives

\begin{equation}
w_{(a)}^{i}(y)z_{b}^{(a)}(x)T_{cd}^{b}(x)w_{(e)}^{c}(x)z_{j}^{(e)}(y)w_{(f)}^{d}(x)z_{k}^{(f)}(y)=T_{jk}^{i}(y)
\end{equation}%
Seperating the variables in (8), (8) is equivalent to

\begin{equation}
z_{a}^{(i)}(x)T_{bc}^{a}(x)w_{(j)}^{b}(x)w_{(k)}^{c}(x)=z_{a}^{(i)}(y)T_{bc}^{a}(y)w_{(j)}^{b}(y)w_{(k)}^{c}(y)
\end{equation}

From (9) we conclude that the expression on the left hand side of (9) is
independent of the variable $x.$ Therefore, this expression is constant on $%
M $ and there exist some constants $C_{(j)(k)}^{(i)}$ such that

\begin{equation}
z_{a}^{(i)}(x)T_{bc}^{a}(x)w_{(j)}^{b}(x)w_{(k)}^{c}(x)=C_{(j)(k)}^{(i)}
\end{equation}
\ 

It is desirable to the express the important formula (10) in a coordinate
free language. Given the parallelizable manifold $(M,w),$ we choose some $%
p\in M$ and define a bilinear form $T(p):T_{p}(M)\times T_{p}(M)\rightarrow
T_{p}(M)$ by the formula

\begin{equation}
\left[ T(p)(\xi ,\mu )\right] ^{i}\overset{def}{=}T_{ab}^{i}(\overline{p}%
)\xi ^{a}\mu ^{b}
\end{equation}%
where $\overline{p}=(\overline{p}^{i})$ are the coordinates of $p$ in the
coordinate system used in (11). Clearly $T(p)$ is skewsymmetric by (4) but
does not necessarily satisfy the Jacobi identity, i.e., it does not define a
Lie algebra. Now $z(p)=(z_{j}^{(i)}(\overline{p}))$ defines an isomorphism $%
\natural z_{p}:T_{p}(M)\rightarrow T_{o}(\mathbb{R}^{n})$ by the formula $%
\natural z_{p}(\xi )\overset{def}{=}z_{a}^{(i)}(\overline{p})\xi ^{a}.$ This
isomorphism extends to the tensor algebra $\natural
z_{p}:T_{p}^{r,s}(M)\rightarrow T_{p}^{r,s}(M).$ Let $A(o)$ be the image of $%
T(p)$ under this isomorphism. Doing this construction at some $q\in M,$ let $%
B(o)$ be the image of $T(q).$ Now we need \textit{not} have $A(o)=B(o)$
because the isomorphism $\natural \varepsilon (p,q):T_{p}(M)\rightarrow
T_{q}(M)$ (which extends to the tensor algebra in the same way) need not map 
$T(p)$ to $T(q).$ Indeed, we have $\natural z_{q}\circ \natural \varepsilon
(p,q)=\natural z_{p}$ by (2) and therefore $\natural z_{q}\circ \natural
\varepsilon (p,q)T_{p}(M)=\natural z_{p}T_{p}(M).$ Thus $\natural
\varepsilon (p,q)T_{p}(M)=T_{q}(M)$ implies $A(o)=B(o).$ However, $\mathfrak{%
R}=0$ implies $\varepsilon $-invariance of $T$ and therefore the
independence of $\natural z_{p}(T(p))=C\in T_{o}^{1,2}(\mathbb{R}^{n})$ on $%
p.$ Now we can easily show that (11) becomes a Lie algebra whose structure
constants with respect to the coordinate basis at $p$ used in (11) are the
constants $C_{(j)(k)}^{(i)}$ defined by (10). Furthermore, (11) is the
localization at $p$ of some $\varepsilon $-invariant Lie algebra of vector
fields on $M$ under the usual bracket.

This coordinate free argument indicates a very simple general principle
underlying the above constancy condition.

\section{A general principle on a LLG}

We observe that $w_{(j)}^{i},$ $1\leq j\leq n=\dim M,$ define $n$
independent vector fields $w_{(j)}$ on $M$ and the tensor $T=(T_{jk}^{i})$
is defined in terms of the first order derivatives of these vector fields,
i.e., $T$ is a first order joint tensorial differential invariant (tdi) of $%
w_{(j)}.$ Now given a smooth manifold $M,$ let $\alpha $ be a tensor field
on $M$ and $\beta $ another tensor field defined using the $k$'th order
derivatives of $\alpha ,$ i.e., $\beta $ is a tdi of $\alpha $ of order $k.$
Here are three well known examples.

1) If $\alpha $ is a $m$-form, its exterior derivative $d\alpha =\beta $ is
a tdi of $\alpha $ of order one.

2) If $\alpha $ is a Riemannian metric, its Riemann curvature tensor, Ricci
tensor and scalar curvature are tdi's of $\alpha $ of order two.

3) If $\alpha $ is an almost complex structure, then its Nijenhuis tensor $%
\beta =N(\alpha )$ is a tdi of order one.

In fact, we can replace $\alpha $ with $r$ tensor fields $\alpha
_{1},...,\alpha _{r}$ and assume that $\beta $ is a joint tdi of $\alpha
_{i} $'s as in our main example of $(M,w).$ Our principle applies also in
this case, but for simplicity, we will assume $r=1$ henceforth.

Now we assume that $\alpha $ is an $\varepsilon $-invariant $(r,s)$-tensor
field on some parallelizable $(M,w).$ By the above argument $\natural
z_{p}(\alpha (p))\in T_{o}^{r,s}(\mathbb{R}^{n})$ does not depend on $p\in
M. $ We call the components of this tensor the structure constants of $%
\alpha .$ As in (10), we can express these components in local coordinates
around $p.$ In this setting, our general principle is

\bigskip

\textbf{GP: }If $\mathfrak{R}=0,$ then any tdi $\beta $ of $\alpha $ of any
order is also $\varepsilon $-invariant.

\bigskip

Indeed, $\mathfrak{R}=0$ implies that the groupoid $\Upsilon $ of $1$-arrows
integrates to the pseudogroup $\mathcal{G}$ whose local diffeomorphisms
leave $\alpha $ invariant since $\alpha $ is $\varepsilon $-invariant.
Differentiation of this invariance condition up to order $k$ shows that $%
\mathcal{G}$ leaves $k$-jets of $\alpha $ invariant for any $k.$ It follows
that $\mathcal{G}$ leaves also $\beta $ invariant. However $1$-arrows
induced by $\mathcal{G}$ belong to $\Upsilon $ which means that $\beta $ is $%
\varepsilon $-invariant. Consequently, the tensor $\natural z_{p}(\beta (p))$
is constant whose components we call the structure constants of $\beta .$

So far so good, but now comes a nontrivial fact.

\begin{proposition}
Let $(M,w,\mathcal{G})$ be a LLG, i.e., some parallelizable $(M,w)$ with $%
\mathfrak{R}=0.$ Let $\alpha $ be an $\varepsilon $-invariant tensor field
and $\beta $ a tdi of $\alpha $ (which is also $\varepsilon $-invariant by
GP). Then the structure constants of $\beta $ can be expressed in terms of
the structure constants of $\alpha $ and the structure constants of the Lie
algebra defined by (11).
\end{proposition}

Proposition 1 holds trivially for $a=w$ and $\beta =T$ by (10). Rather than
attempting a rigorous statement and proof of Proposition 1, we will
illustrate it with three well known examples.

\section{Three examples}

1) Let $\widehat{J}$ be the $2n\times 2n$ matrix whose diagonal consists of $%
n$ blocks of $\widehat{J}_{2}=$ $\left[ 
\begin{array}{cc}
0 & 1 \\ 
-1 & 0%
\end{array}%
\right] .$ Now $\widehat{J}$ defines a linear map of $\mathbb{R}^{2n}$ with
the property $\widehat{J}^{2}=-I_{2n}.$ Given some parallelizable $(M,w),$
we can carry $\widehat{J}$ from $\mathbb{R}^{2n}$ onto $M$ by defining

\begin{equation}
J_{j}^{i}(x)\overset{def}{=}w_{(a)}^{i}(x)\widehat{J}%
_{(b)}^{(a)}z_{j}^{(b)}(x)
\end{equation}%
\ 

Clearly $J=(J_{j}^{i})$ is a 1-1 tensor field on $M,$ that is, a linear map
at the tangent space $T_{p}(M)$ for all $p\in M$ and we easily check that $%
J^{2}=-I,$ i.e., $J$ is an almost complex structure on $M.$ Clearly $%
\widetilde{\nabla }J=0$ or equivalently $J$ is $\varepsilon $-invariant. The
structure constants of $\ J$ are the components of $\widehat{J}.$ We call $J$
the canonical almost complex structure of $(M,w).$

We recall that the Nijenhuis tensor $N(J)$ of $J$ is defined by

\begin{equation}
N(J)_{jk}^{i}=\left[ J_{j}^{a}\frac{\partial J_{k}^{i}}{\partial x^{a}}%
+J_{a}^{i}\frac{\partial J_{j}^{a}}{\partial x^{k}}\right] _{[jk]}=J_{j}^{a}%
\frac{\partial J_{k}^{i}}{\partial x^{a}}+J_{a}^{i}\frac{\partial J_{j}^{a}}{%
\partial x^{k}}-J_{k}^{a}\frac{\partial J_{j}^{i}}{\partial x^{a}}-J_{a}^{i}%
\frac{\partial J_{k}^{a}}{\partial x^{j}}
\end{equation}

It is easy to deduce the expression in (13) in coordinates: Let $\overline{x}%
\in (U,x^{i}),$ $J=(J_{j}^{i}(x))$ be a matrix defined on $(U,x^{i})$
satisfying $J^{2}=-I$ on $U.$ We can always find a coordinate change $%
(x)\rightarrow (y)$ such that $J(\overline{y})=\widehat{J},$ that is, we can
always normalize the value $J(\overline{x}).$ Can we also normalize its
derivative as zero, i.e., can we find a coordinate change $(x)\rightarrow
(y) $ with the property $\frac{\partial J_{j}^{i}(\overline{y})}{\partial
y^{k}}=0?$ Some elementary computations now show that a necessary condition
is $N(J)=0$ at $\overline{x}.$ The local converse of this statement is the
Newlander-Nirenberg theorem: If $N(J)=0$ on $(U,x^{i}),$ then there exists a
coordinate change $(x)\rightarrow (y)$ such that $J(y)=\widehat{J}$ on $%
(V,y^{i}).$

Now we substitute (12) into (13). After some straightforward computation, we
find

\begin{equation}
N(J)_{jk}^{i}=[X_{jk}^{i}]_{[jk]}=X_{jk}^{i}-X_{kj}^{i}
\end{equation}%
where $X_{jk}^{i}$ is the expression

\begin{eqnarray}
&&X_{jk}^{i}\overset{def}{=}w_{(b)}^{a}\widehat{J}_{(c)}^{(b)}z_{j}^{(c)}%
\frac{\partial w_{(d)}^{i}}{\partial x^{a}}\widehat{J}%
_{(e)}^{(d)}z_{k}^{(e)}+w_{(b)}^{a}\widehat{J}%
_{(c)}^{(b)}z_{j}^{(c)}w_{(d)}^{i}\widehat{J}_{(e)}^{(d)}\frac{\partial
z_{k}^{(e)}}{\partial x^{a}}  \notag \\
&&+w_{(b)}^{i}\widehat{J}_{(c)}^{(b)}z_{a}^{(c)}\frac{\partial w_{(d)}^{a}}{%
\partial x^{k}}\widehat{J}_{(e)}^{(d)}z_{j}^{(e)}+w_{(b)}^{i}\widehat{J}%
_{(c)}^{(b)}z_{a}^{(c)}w_{(d)}^{a}\widehat{J}_{(e)}^{(d)}\frac{\partial
z_{j}^{(e)}}{\partial x^{k}}
\end{eqnarray}%
We handle the first term on the right hand side of (15) as follows:

\begin{eqnarray}
w_{(b)}^{a}\widehat{J}_{(c)}^{(b)}z_{j}^{(c)}\frac{\partial w_{(d)}^{i}}{%
\partial x^{a}}\widehat{J}_{(e)}^{(d)}z_{k}^{(e)} &=&w_{(b)}^{a}\widehat{J}%
_{(c)}^{(b)}z_{j}^{(c)}\frac{\partial w_{(d)}^{i}}{\partial x^{a}}\delta
_{(s)}^{(d)}\widehat{J}_{(e)}^{(s)}z_{k}^{(e)}  \notag \\
&=&w_{(b)}^{a}\widehat{J}_{(c)}^{(b)}z_{j}^{(c)}\frac{\partial w_{(d)}^{i}}{%
\partial x^{a}}z_{m}^{(d)}w_{(s)}^{m}\widehat{J}_{(e)}^{(s)}z_{k}^{(e)} 
\notag \\
&=&w_{(b)}^{a}\widehat{J}_{(c)}^{(b)}z_{j}^{(c)}\Gamma _{ma}^{i}w_{(s)}^{m}%
\widehat{J}_{(e)}^{(s)}z_{k}^{(e)} \\
&=&J_{j}^{a}\Gamma _{ba}^{i}J_{k}^{b}  \notag
\end{eqnarray}%
where we substituted from (1) and (3). Similar computation for the second
and third terms give

\begin{equation}
w_{(b)}^{a}\widehat{J}_{(c)}^{(b)}z_{j}^{(c)}w_{(d)}^{i}\widehat{J}%
_{(e)}^{(d)}\frac{\partial z_{k}^{(e)}}{\partial x^{a}}=-J_{j}^{a}J_{b}^{i}%
\Gamma _{ka}^{b}
\end{equation}

\begin{equation}
w_{(b)}^{i}\widehat{J}_{(c)}^{(b)}z_{a}^{(c)}\frac{\partial w_{(d)}^{a}}{%
\partial x^{k}}\widehat{J}_{(e)}^{(d)}z_{j}^{(e)}=J_{a}^{i}\Gamma
_{bk}^{a}J_{j}^{b}
\end{equation}%
For the fourth term we have

\begin{eqnarray}
w_{(b)}^{i}\widehat{J}_{(c)}^{(b)}z_{a}^{(c)}w_{(d)}^{a}\widehat{J}%
_{(e)}^{(d)}\frac{\partial z_{j}^{(e)}}{\partial x^{k}} &=&w_{(b)}^{i}%
\widehat{J}_{(c)}^{(b)}\delta _{(d)}^{(c)}\widehat{J}_{(e)}^{(d)}\frac{%
\partial z_{j}^{(e)}}{\partial x^{k}}  \notag \\
&=&w_{(b)}^{i}\widehat{J}_{(c)}^{(b)}\widehat{J}_{(e)}^{(c)}\frac{\partial
z_{j}^{(e)}}{\partial x^{k}}  \notag \\
&=&-w_{(b)}^{i}\delta _{(e)}^{(b)}\frac{\partial z_{j}^{(e)}}{\partial x^{k}}
\\
&=&-w_{(a)}^{i}\frac{\partial z_{j}^{(a)}}{\partial x^{k}}  \notag \\
&=&-\Gamma _{jk}^{i}  \notag
\end{eqnarray}%
From (15), (16), (17), (18) and (19) we conclude

\begin{equation}
X_{jk}^{i}=J_{k}^{b}J_{j}^{a}\Gamma _{ba}^{i}-J_{j}^{a}J_{b}^{i}\Gamma
_{ka}^{b}+J_{a}^{i}\Gamma _{bk}^{a}J_{j}^{b}-\Gamma _{jk}^{i}
\end{equation}%
and therefore

\begin{eqnarray}
N(J)_{jk}^{i} &=&X_{jk}^{i}-X_{kj}^{i} \\
&=&T_{ab}^{i}J_{k}^{a}J_{j}^{b}+J_{a}^{i}T_{bk}^{a}J_{j}^{b}-J_{a}^{i}T_{bj}^{a}J_{k}^{b}-T_{jk}^{i}
\notag
\end{eqnarray}

We observe that (21) holds on any parallelizable $(M,w).$ If $\mathfrak{R}%
=0, $ note that $\varepsilon $-invariance of $N(J),$ a consequence of GP,
can be seen directly from (21) by applying $\widetilde{\nabla }$ to both
sides of (21) in view of (5). Using (11), we can write (21) in a coordinate
free form as

\begin{equation}
N(J)(u,v)=[J(u),J(v)]+J[J(u),v]-J[J(v),u]-[u,v]
\end{equation}%
Note that $[$ $,$ $]$ is the usual bracket of vector fields in (22) which
localizes at any point according to (11). Generalizing this bracket to Fr%
\"{o}hlicher-Nijenhuis bracket, we recall that (22) can be taken as the
modern definition of the Nijenhuis tensor on any almost complex manifold. We
observe that $N(J)_{ak}^{a}=-2T_{ak}^{a}$ by (21) so that $N(J)$ defines the
same secondary characteristic classes as $T$ if $\mathfrak{R}=0$ (see
Chapter 12 of [O1]).

Now if $\mathfrak{R}=0,$ then $N(J)$ is $\varepsilon $-invariant by GP and
therefore $\natural z_{p}\left( N(J)(p)\right) \in T_{o}^{1,2}(\mathbb{R}%
^{2n})$ is constant. In coordinates, this means

\begin{equation}
z_{a}^{(i)}N(J)_{bc}^{a}w_{(j)}^{b}w_{(k)}^{c}=N(J)_{(j)(k)}^{(i)}
\end{equation}%
in the same way as (10). Evaluating (21) at $p$ and translating (21) to $%
o\in \mathbb{R}^{2n}$ by the isomorphism $\natural z_{p}$ amounts to
contracting both sides of (21) with the components of $z$ and $w.$ The
result is

\begin{equation}
N(J)_{(j)(k)}^{(i)}=C_{(a)(b)}^{(i)}\widehat{J}_{(k)}^{(a)}\widehat{J}%
_{(j)}^{(b)}+C_{(b)(k)}^{(a)}\widehat{J}_{(a)}^{(i)}\widehat{J}%
_{(j)}^{(b)}-C_{(b)(j)}^{(a)}\widehat{J}_{(a)}^{(i)}\widehat{J}%
_{(k)}^{(b)}-C_{(j)(k)}^{(i)}
\end{equation}%
which expresses the structure constants of $N(J)$ in terms of the structure
constants of $J$ and the Lie algebra (11). We observe that (21) and (24) are
essentially the same formulas but (21) holds on $(M,w)$ whereas (24) needs $%
\mathfrak{R}=0.$

2) Let $\widehat{\varpi }$ be the $2n\times 2n$ matrix $\left[ 
\begin{array}{cc}
0 & I_{n} \\ 
-I_{n} & 0%
\end{array}%
\right] .$ We define the skewsymmetric and nondegenerate 2-form $\varpi
=(\varpi _{ij})$ on $(M,w)$ by

\begin{equation}
\widehat{\varpi }_{(a)(b)}z_{i}^{(a)}z_{j}^{(b)}=\varpi _{ij}(x)
\end{equation}

We call $\varpi $ the canonical almost symplectic structure of $(M,w)$ which
is $\varepsilon $-invariant with structure constants as the components of $%
\widehat{\varpi }.$ An easy computation shows that a necessary condition to
normalize the $1$-jet of $\varpi $ as zero is the vanishing of the exterior
derivative

\begin{equation}
\left( d\varpi \right) _{kij}=\frac{\partial \varpi _{ij}}{\partial x^{k}}-%
\frac{\partial \varpi _{kj}}{\partial x^{i}}-\frac{\partial \varpi _{ik}}{%
\partial x^{j}}
\end{equation}%
(omitting the factor $\frac{1}{3})$ at that point and the local converse of
this statement is the Darboux theorem. We substitute (25) into (26) and what
we find is

\begin{equation}
(d\varpi )_{kij}=T_{ki}^{a}\varpi _{ja}-T_{kj}^{a}\varpi
_{ia}-T_{ji}^{a}\varpi _{ka}
\end{equation}%
or

\begin{equation}
d\varpi (\xi ,\gamma ,\rho )=\varpi (\rho ,[\xi ,\gamma ])-\varpi (\gamma
,[\xi ,\rho ])-\varpi (\xi ,[\rho ,\gamma ])
\end{equation}%
which holds on $(M,w).$ Note that (28) is "half" of the well known formula
for the exterior derivative (the other half vanishes since $\varpi $ is $%
\widetilde{\nabla }$-parallel and $\widetilde{\nabla }$ is "formal Lie
derivative" ([O2]). If $\mathfrak{R}=0,$ (27) gives

\begin{equation}
(d\varpi )_{(i)(j)(k)}=T_{(k)(i)}^{(a)}\widehat{\varpi }%
_{(j)(a)}-T_{(k)(j)}^{(a)}\widehat{\varpi }_{(i)(a)}-T_{(j)(i)}^{(a)}%
\widehat{\varpi }_{(k)(a)}
\end{equation}

3) This is by far the most important example and is studied to some extent
in [O3] (we warn the reader that some formulas in [O3] are not correct as
they stand). The local computations are quite more involved compared to the
above examples since Rieman curvature tensor is a \textit{second order }tdi
of the canonical metric defined by

\begin{equation}
g_{ij}(x)\overset{def}{=}\widehat{g}_{(a)(b)}z_{i}^{(a)}z_{j}^{(b)}
\end{equation}%
where $\widehat{g}_{ij}\overset{def}{=}\delta _{ij}$ and therefore $%
g_{ij}=z_{i}^{(a)}z_{j}^{(a)}$ with a summation over $a.$ We will leave the
details of these computations to the interested reader. However, there is an
important issue here which has been lingering for a long time and needs to
be settled: Since $g$ is $\varepsilon $-invariant, the scalar curvature $c(x)
$ is also $\varepsilon $-invariant by GP if we assume $\mathfrak{R}=0.$
However a function is $\varepsilon $-invariant if and only if it is constant
and there the $c(x)$ is constant. On the other hand, $c$ can be expressed in
terms of the structure constants of the Lie algebra defined by (11)
according to Proposition 1. We invite the interested reader to derive this
expression. Therefore, Pommaret's reasoning that "a single constant can not
be related to any Lie algebra" which he uses to support his conceptual
framework in many of his works (see, for instance [P], pg.29), is not
justified. He also claims in various places (see again [P], pg.29) that the
modern concept of curvature must be revisited due to a confusion between the
Janet and Spencer sequences. We believe that there are two concepts of
curvatures in geometry. The first one is well known and is widely accepted
by the geometers: It originates from Riemannian geometry where the central
concept is a connection and curvature is an invariant of a connection
together with a very intriguing invariant called torsion. As a peculiarity
among others, only connections on "certain" principal bundles have torsion!
The second one emerges from Lie theory as a "local obstruction to the
integration of the groupoid of $k$-arrows to a pseudogroup" or equivalently,
as a local obstruction to the formal integrability of certain geometric
PDE's where a "connection" appears as part of the definition of these PDE's.
This second curvature is formulated, in its simplest incarnation $k=0,$ in
Part I of [O1] and studied further in [O1], [O2] where torsion becomes the
measure of the difference between the "left-right actions". Even though the
distinction between these two curvatures may seem to be only a matter of
interpretation at first sight, it has some quite unexpected consequences.

Given some parallelizable $(M,w)$ we can now use the notation $(M,w,J,\varpi
,g)$ where the last three structures are canonically defined as above and
any two of them determines the other. This raises a natural question: What
are the relations between the tdi's of these structures? Another question is
how these arguments generalize to prehomogeneous geometries in [O1], [O2].

\bigskip

We will conclude this note by expressing our conviction that a
parallelizable manifold and more generally a prehomogeneous geometry is an
extremely rich geometric structure encompassing most (if not all) of
differential geometry originating from Klein's Erlangen Program.

\bigskip

\textbf{References}

\bigskip

[O1] E.H.Orta\c{c}gil: An Alternative Approach to Lie Groups and Geometric
Structures, OUP, 2018

[O2] E.H.Orta\c{c}gil: Curvature without connection, arXiv: 2003.06593, 2020

[O3] E.H.Orta\c{c}gil: The canonical geometry of a local Lie group, arXiv:
2004.04029, 2020

[P] J.F.Pommaret: Partial Differential Equations and Group Theory, New
Perspectives for Applications, Kluwer Academic Publishers, 1994

\bigskip

Erc\"{u}ment Orta\c{c}gil

ortacgile@gmail.com

\end{document}